\documentclass[final,a4paper]{article}
\usepackage{fullpage}
\usepackage[latin1]{inputenc}
\usepackage[dvips]{graphicx}
\usepackage{amsmath,amsfonts,amsthm}
\usepackage{xspace}
\newcommand{\conjugate}[1]{\overline{#1}}
\newcommand{\involute}[2]{\overline{#1\,}^{#2}}
\renewcommand{\v}[1]{\boldsymbol{#1}} 
\renewcommand{\i}{\v{i}}
\renewcommand{\j}{\v{j}}
\renewcommand{\k}{\v{k}}
\newcommand{\vmu}{\v\mu}
\newcommand{\vnu}{\v\nu}
\newcommand{\ie}{\textit{i.e.}}
\renewcommand{\H}{\ensuremath{\mathbb{H}}\xspace}
\newtheorem{lemma}{Lemma}
\newtheorem{theorem}{Theorem}
\newtheorem{axiom}{Axiom}
\newtheorem{corollary}{Corollary}
\title{Quaternion Involutions}
\author{Todd A. Ell\thanks{Dr T. A. Ell is a Visiting Fellow at the University of Essex, funded by grant
                           number GR/S58621 from the United Kingdom Engineering and Physical Sciences
                           Research Council.}\\
        5620 Oak View Court, Savage, MN 55378-4695, USA.\\
        Email:~T.Ell@IEEE.org.\\
\and Stephen J. Sangwine\\
     Department of Electronic Systems Engineering,\\
     University of Essex, Wivenhoe Park,\\
     Colchester, CO7 9EU, United Kingdom.\\
     Email:~S.Sangwine@IEEE.org.}
\begin{document}
\maketitle
\begin{abstract}
An involution is usually defined as a mapping that is its own inverse. In this paper, we
study quaternion involutions that have the additional properties of distribution over
addition and multiplication. We review formal axioms for such involutions, and
we show that the quaternions
have an infinite number of involutions. We show that the conjugate of a quaternion may be expressed
using three mutually perpendicular involutions. We also show that any set of three
mutually perpendicular quaternion involutions is closed under composition. Finally, we
show that projection of a vector or quaternion can be expressed concisely using
involutions.
\end{abstract}

\section{Introduction}
Involutions are usually defined simply as self-inverse mappings. A trivial example
is conjugation of a complex number, which is obviously self-inverse.
In this paper we consider involutions of the quaternions, that
is functions of a quaternion variable that are self-inverse. Quaternion conjugation
is an obvious involution, but it is not the only quaternion involution. In fact,
the quaternions have an infinite number of involutions, as we show. The paper begins
by reviewing the classical basics of quaternions, and then presents axioms for
involutions which go beyond the simple definition of a self-inverse mapping. Section
\ref{qinvolutions} then presents the quaternion involutions, and section \ref{qinvolutionprops}
presents their properties. Section \ref{conjugate} discusses the quaternion conjugate
and shows that it may be expressed using three mutually perpendicular quaternion
involutions. Finally, section \ref{projection} shows that the projection of a vector or
quaternion may be expressed using involutions.

\section{Basics of quaternions}
A quaternion may be represented in Cartesian form $q = w + \i x + \j y + \k z$
where $\i$, $\j$ and $\k$ are mutually perpendicular unit vectors obeying the
multiplication rules below discovered by Hamilton in 1843 \cite{Hamilton:1866},
and $w$, $x$, $y$, $z$, are real.
\begin{gather}
\label{hrules} \i^2 = \j^2 = \k^2 = \i\j\k = -1\
\end{gather}
The conjugate of a quaternion is given by $\conjugate{q} = w - \i x - \j y - \k z$.

The quaternion algebra \H is a normed division algebra. The modulus of a quaternion is
the square root of its norm: $|q| = \sqrt{w^2 + x^2 + y^2 + z^2}$, and every non-zero
quaternion has a multiplicative inverse given by its conjugate divided by its norm:
$q^{-1} = \conjugate{q}/|q|^2 = (w - \i x - \j y - \k z)/(w^2 + x^2 + y^2 + z^2)$.
For a more detailed exposition of the basics of quaternions, we refer the reader
to Coxeter's 1946 paper \cite{Coxeter:1946}.

An alternative and much more powerful representation for a quaternion is as a combination
of a \emph{scalar} and a \emph{vector} part, analogous to a complex number, and this
representation will be employed in the rest of the paper: $q = a + \vmu b$, where $\vmu$
is a unit vector, and $a$ and $b$ are real. $b$ is the modulus of the vector part of the
quaternion and $\vmu$ is its direction. In terms of the Cartesian representation:
\begin{equation}
\label{svform}
a = w,\qquad b = \sqrt{x^2 + y^2 + z^2},\qquad \vmu = \frac{\i x + \j y + \k z}{b}
\end{equation}
\begin{lemma}
\label{unitvector}
The square of any unit vector is $-1$.
\end{lemma}
\begin{proof}
Let $\vmu$ be an arbitrary unit vector as defined in equation \ref{svform}.
Its square is given by:
\begin{equation*}
\begin{split}
\vmu^2 &= \frac{(\i x + \j y + \k z)^2}{x^2 + y^2 + z^2}\\
       &= \frac{\i^2 x^2 + \j^2 y^2 + \k^2 z^2 +
                \i\j x y + \j\i x y + \i\k x z + \k\i x z + \j\k y z + \k\j y z}{x^2 + y^2 + z^2}\\
\end{split}
\end{equation*}
Applying the rules in equation \ref{hrules} we get:
$\vmu^2 = \frac{-x^2 -y^2 -z^2}{x^2 + y^2 + z^2} = -1$
\end{proof}
A corollary of Lemma \ref{unitvector} is that there are an infinite number of
solutions to the equation $x^2 = -1$.

The conjugate of a quaternion in complex form is $\conjugate{q} = a - \vmu b$.
Geometrically, this is obviously a reversal of the direction
of the vector part. The quaternion conjugate has properties
analogous to those of the complex conjugate, with one minor
exception: the quaternion conjugate is an anti-involution whereas
the complex conjugate is an involution. We define these terms in
the next section, and we return to this point with Theorem \ref{anticonjugate}
in section \ref{conjugate}. The product of a quaternion with its conjugate
gives the norm, or square of the modulus. This follows directly
from Lemma \ref{unitvector}:
$(a + \vmu b)(a - \vmu b) = a^2 - \vmu^2 b^2 = a^2 + b^2 = w^2 + x^2 + y^2 + z^2$.

\section{Involutions}
\label{involutions}
The formal definition of an involution is not easy to find (most mathematical
reference works define it simply as a mapping which is its own inverse) but
\cite{KluwerEncylMaths:1990} gives a reasonably authoritative statement from
which we reproduce the following axioms. It is clear from what follows in the
paper that all three of these axioms are important, otherwise it is possible
to define trivial self-inverse mappings, which have uninteresting properties.
%
%
\newcommand{\involution}[1]{f(#1)}
We denote an arbitrary involution by the mapping $x \rightarrow \involution{x}$.
\begin{axiom}
  \label{inverse}
  $\involution{\involution{x}} = x$. An involution is its own inverse.
\end{axiom}
\begin{axiom}
  \label{linearity}
  An involution is linear:
  $\involution{x_1 + x_2} = \involution{x_1} + \involution{x_2}$ and
  $\lambda\involution{x} = \involution{\lambda x}$ where $\lambda$ is
  a real constant.
\end{axiom}
\begin{axiom}
  \label{product}
  $\involution{x_1 x_2} = \involution{x_1}\involution{x_2}$. If the terms
  on the right must be reversed (which can only be necessary if $x_1$ and $x_2$
  do not commute), then $\involution{\ }$ is an \emph{anti-involution}.
\end{axiom}
\section{Quaternion involutions}
\label{qinvolutions}

Involutions over the quaternion field have been published by Chernov \cite{Chernov:1995},
and used by Bülow \cite{Bulow:1999,BulowSommer:2001}, but with an important difference
from the involutions given in this paper: we show here that there are an
infinite number of involutions over the quaternion field
whereas Chernov and Bülow wrote of
only three (plus conjugation). Chernov defined the
following involutions, which we generalize in Theorem \ref{general}:
\begin{equation}
\begin{split}
\alpha(q) &= - \i q \i = w + \i x - \j y - \k z\\
\beta(q)  &= - \j q \j = w - \i x + \j y - \k z\\
\gamma(q) &= - \k q \k = w - \i x - \j y + \k z\\
\end{split}
\end{equation}
They also showed that the quaternion conjugate can be expressed in terms of these three
involutions, a result that we generalize in Theorem \ref{genconj}.

\begin{theorem}
\label{general}
The mapping $q \rightarrow -\vnu q \vnu$ where $q$ is an arbitrary quaternion is
an involution for any unit vector $\vnu$.
\end{theorem}
\begin{proof}
Axiom \ref{inverse} is easily shown to be satisfied using Lemma \ref{unitvector}:
$-\vnu (-\vnu q \vnu) \vnu = (-1)q(-1) = q$. Axiom \ref{linearity} is seen to be
satisfied from: $-\vnu (q_1 + q_2) \vnu = -\vnu q_1 \vnu - \vnu q_2 \vnu$ (multiplication
of quaternions is distributive over addition). Since reals commute with quaternions, the
second part of the axiom is trivially seen.

Axiom \ref{product} can be shown to be satisfied as follows:
\begin{equation*}
\begin{split}
\involution{q_1}\involution{q_2} & = (-\vnu q_1 \vnu)(-\vnu q_2 \vnu)\\
                                 & = \vnu q_1 \vnu \vnu q_2 \vnu\\
                                 & = \vnu q_1 (-1) q_2 \vnu\\
                                 & = -\vnu q_1 q_2 \vnu = \involution{q_1 q_2}\\
\end{split}
\end{equation*}
\end{proof}

Note that a mapping $q \rightarrow \vnu_1 q \vnu_2$ ($\vnu_1 \neq \vnu_2$) is its
own inverse, but is not an involution as considered in this paper, because it does
not satisfy axiom \ref{product}.

In what follows we introduce a new notation for involutions using an overbar with a
subscript unit vector, thus $\involute{q}{\vnu} = -\vnu q \vnu$. We refer to the direction
defined by $\vnu$ as the \emph{axis of involution}. The fact that we use an overbar to
denote involutions as well as the conjugate is not coincidental and we shall see that
there is a close relationship between involutions and conjugation.

\section{Properties of quaternion involutions}
\label{qinvolutionprops}
\begin{lemma}
\label{unitvectorproduct}
The product of any two vectors with arbitrary non-zero norms is a quaternion. The scalar
part is minus the inner or scalar product of the two vectors and the vector part is the
vector product of the two vectors. Reversing the order of the product conjugates the
resulting quaternion.
\end{lemma}
\begin{proof}
Let $\vmu_1 = \i x_1 + \j y_1 + \k z_1$ and $\vmu_2 = \i x_2 + \j y_2 + \k z_2$.
Their product is given by:
\begin{equation*}
\begin{split}
\vmu_1 \vmu_2 & = (\i x_1 + \j y_1 + \k z_1)(\i x_2 + \j y_2 + \k z_2)\\
              & = \i^2 x_1 x_2 + \j^2 y_1 y_2 + \k^2 z_1 z_2\\
              & \quad + \i\j x_1 y_2 + \j\i y_1 x_2 + \i\k x_1 z_2 + \k\i z_1 x_2 + \j\k y_1 z_2 + \k\j z_1 y_2\\
              & = - (x_1 x_2 + y_1 y_2 + z_1 z_2)\\
              & \quad + \i(y_1 z_2 - z_1 y_2) + \j(z_1 x_2 - x_1 z_2) + \k(x_1 y_2 - y_1 x_2)
\end{split}
\end{equation*}
Changing the order of the product changes the order of all the products of two unit
vectors $\i$, $\j$ and $\k$. Since $\i\j = -\j\i$ and so on, this negates all the
components of the vector part of the result. The scalar part, $-(x_1 x_2 + y_1 y_2 + z_1 z_2)$,
is unchanged. Thus reversing the order of the product conjugates the result, as stated.
The scalar and vector parts can be seen to be equal to minus the scalar product and
the vector product respectively, according to standard definitions of these products.
\end{proof}

\begin{lemma}
\label{perpvec}
The product of two perpendicular vectors changes sign if the order of the product
is reversed.
\end{lemma}
\begin{proof}
Lemma \ref{unitvectorproduct} identified the scalar part of the result
$-(x_1 x_2 + y_1 y_2 + z_1 z_2)$ as minus the inner product of the two
vectors. Since the inner product is zero in the case of perpendicular vectors, the
product of perpendicular vectors is a vector and it is the vector product of the
two vectors. It follows that this vector changes sign (reverses) if the
order of the product is reversed.
\end{proof}
\begin{theorem}
\label{composition}
Composition: the composition of two perpendicular involutions is commutative. That is:
$\involute{\involute{q}{\vnu_1}}{\vnu_2} = \involute{\involute{q}{\vnu_2}}{\vnu_1}$ where
$\vnu_1\perp\vnu_2$.
\end{theorem}
\begin{proof}
\[
\involute{\involute{q}{\vnu_1}}{\vnu_2} = -\vnu_2(-\vnu_1 q \vnu_1)\vnu_2 = \vnu_2\vnu_1 q \vnu_1\vnu_2
\]
and by Lemma \ref{perpvec} we can reverse the order of the pairs of unit vectors
if we change their signs:
\[
\involute{\involute{q}{\vnu_1}}{\vnu_2} = (-\vnu_1\vnu_2) q (-\vnu_2\vnu_1) = -\vnu_1(-\vnu_2 q \vnu_2)\vnu_1 = \involute{\involute{q}{\vnu_2}}{\vnu_1}
\]
\end{proof}
\begin{theorem}
\label{perpinvolutions}
Double Composition: given a set of three mutually perpendicular unit vectors,
$\vnu_1$, $\vnu_2$, $\vnu_3$, such that $\vnu_1 \vnu_2 = \vnu_3$, then
\[
\involute{\involute{q}{\vnu_1}}{\vnu_2} = \involute{q}{\vnu_3}
\]
\end{theorem}
\begin{proof}
\begin{equation*}
\begin{split}
\involute{\involute{q}{\vnu_1}}{\vnu_2} & = -\vnu_2 (-\vnu_1 q \vnu_1) \vnu_2 = \vnu_2\vnu_1 q \vnu_1\vnu_2\\
                                        & = -\vnu_1\vnu_2 q \vnu_1\vnu_2 \quad\text{by Lemma \ref{perpvec}}\\
                                        & = -\vnu_3 q \vnu_3 = \involute{q}{\vnu_3}
\end{split}
\end{equation*}
\end{proof}
\begin{corollary}
Triple Composition: the composition of three mutually perpendicular involutions is an identity.
\end{corollary}
\begin{proof}
From Theorem \ref{perpinvolutions}:
$\involute{\involute{q}{\vnu_1}}{\vnu_2} = \involute{q}{\vnu_3}$.
Apply an involution about $\vnu_3$ to both sides:
$\involute{\involute{\involute{q}{\vnu_1}}{\vnu_2}}{\vnu_3} = \involute{\involute{q}{\vnu_3}}{\vnu_3} = q$
\end{proof}
Thus we see that a set of three mutually perpendicular involutions
(that is involutions about a set of three mutually perpendicular
axes) is closed under composition of the involutions
\emph{and by Theorem \ref{composition} the order of the composition is unimportant.}

\begin{figure}
\centerline{\includegraphics[width=0.4\textwidth]{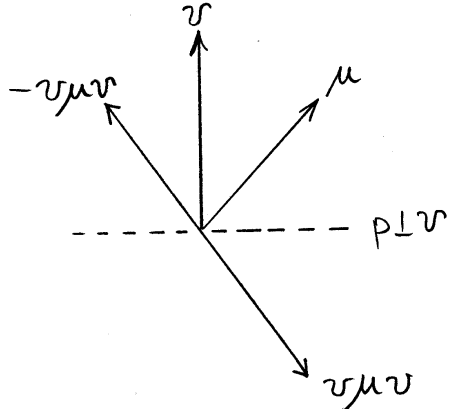}}
\caption{\label{reflect}Reflection of a vector $\vmu$ in a line defined by a unit vector $\vnu$.
                        $p$ is a plane perpendicular to $\vnu$ seen edge on.}
\end{figure}

We now present a geometric interpretation of quaternion involution.
\begin{theorem}
\label{geometry} Given an arbitrary quaternion $q = a + \vmu b$,
an involution $\involute{q}{\vnu}$ leaves the scalar part of $q$
(that is, $a$) invariant, and reflects the vector part of $q$
(that is, $\vmu b$) across the line defined by the
axis of involution $\vnu$. (Equivalently, the vector part of $q$
is rotated by $\pi$ radians about the axis of involution $\vnu$.)
\end{theorem}
\begin{proof}
\begin{equation*}
\begin{split}
\involute{q}{\vnu} & = -\vnu(a + \vmu b)\vnu = -\vnu a\vnu -\vnu\vmu b\vnu = -\vnu^2 a -\vnu\vmu\vnu b
\intertext{and, since $\vnu$ is a unit vector, by Lemma \ref{unitvector}:}
\involute{q}{\vnu} & = a - \vnu\vmu\vnu b
\end{split}
\end{equation*}
We recognise $\vnu\vmu\vnu$ to be a reflection of $\vmu$ in the plane $p$ normal to $\vnu$
as shown by Coxeter \cite[Theorem 3.1]{Coxeter:1946}. Therefore $-\vnu\vmu\vnu$ is a
reflection of $\vmu$ in the line defined by $\vnu$ as shown in Figure \ref{reflect}.
The result of the reflection of the vector part remains a vector, and therefore the
scalar part remains unchanged, as shown, and as stated.
\end{proof}
\begin{corollary}
\label{invparq} An involution applied to a quaternion with vector
part parallel to the involution axis is an identity,
\ie, $\involute{a + \vmu b}{\vnu} = a + \vmu b$, where $\vnu\parallel\vmu$.
\end{corollary}
\begin{proof}
From Theorem \ref{geometry}: $\involute{q}{\vnu} = a - \vnu\vmu\vnu b$. Since both $\vnu$
and $\vmu$ are unit vectors, and are parallel, $\vnu = \pm\vmu$ and the result on the right
reduces to $a + \vmu b = q$ in both cases.
\end{proof}
\begin{corollary}
\label{invperpq} An involution applied to a quaternion with vector
part perpendicular to the involution axis conjugates the quaternion,
\ie, $\involute{a + \vmu b}{\vnu} = a - \vmu b$, where $\vnu\perp\vmu$.
\end{corollary}
\begin{proof}
From Theorem \ref{geometry}: $\involute{q}{\vnu} = a - \vnu\vmu\vnu b$.
If $\vnu\perp\vmu$ we can use Lemma \ref{perpvec} to reverse the order of the two
unit vectors, thus obtaining $a + \vnu^2\vmu b = a - \vmu b = \conjugate{q}$ as stated.
\end{proof}
\begin{lemma}
\label{argvecprod} The product of two unit vectors is a quaternion
with argument equal to the angle between the two vectors.
\end{lemma}
The sign of the argument is significant, because it depends on the
ordering of the two vectors (the angle is measured from the first
vector to the second).
\begin{proof}
Lemma \ref{unitvectorproduct} identified the scalar and vector parts
of the product of two vectors with minus the inner product, and the vector product
respectively of the two vectors. Since these products are given for unit vectors
by $\cos\theta$ and $\vmu\sin\theta$, where $\vmu$ is perpendicular to the plane
containing the two vectors, we may write the product of two unit vectors as:
$-\cos\theta + \vmu\sin\theta = -\exp(-\vmu\theta)$.
\end{proof}
\begin{theorem}
The composition of two involutions is a rotation of the vector part of the quaternion
operated upon about an axis normal to the plane containing the two axes of involution.
The angle of rotation is twice the angle between the two involution axes.
In the case where the two involutions are perpendicular the composite result is of course
an involution, which is a rotation of the vector part of the quaternion by $\pi$.
\end{theorem}
\begin{proof}
Let $\vnu_1$ and $\vnu_2$ be two unit vectors and $q$ be an arbitrary quaternion.
Then the composition of two involutions about $\vnu_1$ and $\vnu_2$ is given by:
\[
\involute{\involute{q}{\vnu_1}}{\vnu_2} = -\vnu_2(-\vnu_1 q \vnu_1)\vnu_2 = \vnu_2\vnu_1 q \vnu_1\vnu_2
\]
From Lemma \ref{unitvectorproduct} we can write the result on the right as $p q
\conjugate{p}$, where $p = \vnu_2\vnu_1$ is a unit quaternion\footnote{$p$ is
a unit quaternion, because it is the product of two unit vectors. This follows from
the fact that the quaternion algebra is a normed algebra.}.
Separating $q = a + \vmu b$ into its scalar and vector parts we obtain:
\[
\involute{\involute{q}{\vnu_1}}{\vnu_2} = p \conjugate{p} a + p \vmu b \conjugate{p} = a + (p\vmu\conjugate{p}) b
\]
We recognise the term $p\vmu\conjugate{p}$ as a rotation as given by Coxeter \cite[Theorem
3.2]{Coxeter:1946}. The axis of rotation is given by the vector part of $p$, and the angle
of rotation is twice the argument of $p$. Therefore, from Lemma \ref{unitvectorproduct} we
know that the axis of rotation is perpendicular to the plane containing the two vectors,
and from Lemma \ref{argvecprod} we know that the angle of rotation is twice the angle
between the two vectors.
\end{proof}

\section{The quaternion conjugate}
\label{conjugate}
\begin{theorem}
\label{anticonjugate}
The quaternion conjugate is an anti-involution.
\end{theorem}
\begin{proof}
The definition of the conjugate of a quaternion $q = a + \vmu b$ is $\conjugate{q} = a -
\vmu b$. We have to show that this satisfies the three axioms given in section \ref{involutions},
and that in Axiom \ref{product} we have to reverse the terms on the
right-hand side.

It is easily seen that the quaternion conjugate satisfies Axiom \ref{inverse}.
To demonstrate that the quaternion conjugate satisfies Axiom \ref{linearity}, let
$q_1 = a + \vmu_1 b$ and $q_2 = c + \vmu_2 d$. Then $q_1 + q_2 = (a + c) + (\vmu_1 b + \vmu_2 d)$.
Since reversing two vectors also reverses their sum, we see that the quaternion conjugate
is distributive over addition, as required. The second part of the axiom is easily seen.

To show that the quaternion conjugate satisfies Axiom \ref{product}, we state
the required equality and then demonstrate that it is satisfied by expanding the left and
right hand sides until identical:
\begin{equation*}
\begin{split}
\conjugate{q_1 q_2}                      & = \conjugate{q_2}\,\conjugate{q_1}\\
\conjugate{(a + \vmu_1 b)(c + \vmu_2 d)} & = (c - \vmu_2 d)(a - \vmu_1 b)\\
\intertext{Expanding the left hand side, we employ Axiom \ref{linearity}:}
a c - \vmu_1\vmu_2 b d - \vmu_1 b c - \vmu_2 a d & = a c + \vmu_2\vmu_1 b d - \vmu_1 b c - \vmu_2 a d\\
\intertext{and using Lemma \ref{perpvec} we change the order of the two vectors
           in the second term on the right-hand side to obtain the required result:}
a c - \vmu_1\vmu_2 b d - \vmu_1 b c - \vmu_2 a d & = a c - \vmu_1\vmu_2 b d - \vmu_1 b c - \vmu_2 a d\\
\end{split}
\end{equation*}
\end{proof}
We now show that the quaternion conjugate can be expressed using the sum of three mutually
perpendicular involutions.
\begin{lemma}
\label{treperpinv}
The sum of three mutually perpendicular involutions applied to a vector negates the
vector (reverses its direction). That is, given a set of three mutually perpendicular
unit vectors as in Theorem \ref{perpinvolutions} and an arbitary vector $\vmu$:
\begin{equation}
\label{treperinveqn}
\involute{\vmu}{\vnu_1} + \involute{\vmu}{\vnu_2} + \involute{\vmu}{\vnu_3} = -\vmu
\end{equation}
\end{lemma}
\begin{proof}
\newcommand{\veta}{\v\eta}
Let $\vmu = \veta_1 + \veta_2 + \veta_3$ where $\veta_i
\parallel \vnu_i, i \in \{1, 2, 3\}$.
In other words,
resolve $\vmu$ into three vectors\footnote{The three vectors $\veta_i$ are not, in general,
of unit modulus.} parallel to the three
mutually perpendicular unit vectors $\vnu_1$, $\vnu_2$ and
$\vnu_3$. Substitute this representation of $\vmu$ into the
left-hand side of Equation \ref{treperinveqn}:
\begin{equation*}
\involute{\veta_1 + \veta_2 + \veta_3}{\vnu_1} +
\involute{\veta_1 + \veta_2 + \veta_3}{\vnu_2} +
\involute{\veta_1 + \veta_2 + \veta_3}{\vnu_3} = -\vmu
\end{equation*}
Axiom \ref{linearity} allows us to apply the involutions separately to the three components:
\begin{equation*}
\involute{\veta_1}{\vnu_1} + \involute{\veta_2}{\vnu_1} + \involute{\veta_3}{\vnu_1} +
\involute{\veta_1}{\vnu_2} + \involute{\veta_2}{\vnu_2} + \involute{\veta_3}{\vnu_2} +
\involute{\veta_1}{\vnu_3} + \involute{\veta_2}{\vnu_3} + \involute{\veta_3}{\vnu_3} = -\vmu
\end{equation*}
We now make use of Corollaries \ref{invparq} and \ref{invperpq}. In this case we are
applying them to a vector, so the first states that an involution with axis parallel to a
vector is an identity, and the second states that an involution with axis perpendicular to
the vector reverses, or negates, the vector:
\begin{equation*}
\veta_1 - \veta_2 - \veta_3 - \veta_1 + \veta_2 - \veta_3 - \veta_1 - \veta_2 + \veta_3 = -\vmu
\end{equation*}
and cancelling out, we obtain: $- \veta_1 - \veta_2 - \veta_3 = -\vmu$, which is the assumption
we made at the start of the proof.
\end{proof}
The following theorem is a generalization of a similar result given in \cite[Definition 2.2, p.12]{Bulow:1999}.
\begin{theorem}
\label{genconj}
Given a set of three mutually perpendicular unit vectors as in Theorem \ref{perpinvolutions},
the conjugate of $q$ may be expressed as:
\begin{equation}
\label{conjinvol}
\conjugate{q} 
                =  \frac{1}{2}\left(\involute{q}{\vnu_1} + \involute{q}{\vnu_2} + \involute{q}{\vnu_3} -q \right)
\end{equation}
\end{theorem}
\begin{proof}
Let $q = a + \vmu b$. Substituting this expression for $q$ into the right-hand side of
Equation \ref{conjinvol} we obtain:
\begin{equation*}
\conjugate{q} = \frac{1}{2}\left(\involute{a + \vmu b}{\vnu_1}
                               + \involute{a + \vmu b}{\vnu_2}
                               + \involute{a + \vmu b}{\vnu_3}
                               -     \left(a + \vmu b\right)
                          \right)
\end{equation*}
We now apply the three involutions separately to the components of $q$ using Axiom
\ref{linearity}, and noting from Theorem \ref{geometry} that the scalar part $a$ is
invariant under involutions:
\begin{equation*}
\conjugate{q} = \frac{1}{2}\left(a + \involute{\vmu}{\vnu_1} b
                               + a + \involute{\vmu}{\vnu_2} b
                               + a + \involute{\vmu}{\vnu_3} b
                               - a -           \vmu b
                          \right)
\end{equation*}
Gathering terms together and factoring out $b$:
\begin{equation*}
\conjugate{q} = a +\frac{1}{2}\left(\involute{\vmu}{\vnu_1}
                                 +  \involute{\vmu}{\vnu_2}
                                 +  \involute{\vmu}{\vnu_3}
                                 -            \vmu  \right) b
\end{equation*}
and the right-hand side is equal to $a - \vmu b$ by Lemma \ref{treperpinv}.
\end{proof}

\section{Projection using involutions}
\label{projection}
Finally, we now demonstrate the utility of quaternion involutions by presenting formulae for
projection of a vector into or perpendicular to a given direction. These results have
been published in \cite{SangwineEll:2000c}, but without explicit use of involutions.
\begin{theorem}
\label{projections}
An arbitrary vector $\vmu$ may be resolved into two components parallel to,
and perpendicular to, a direction in 3-space defined by a unit vector $\vnu$:
\begin{equation*}
\vmu_{\parallel\vnu} = \frac{1}{2}(\vmu + \involute{\vmu}{\vnu})\qquad
\vmu_{\perp\vnu}     = \frac{1}{2}(\vmu - \involute{\vmu}{\vnu})
\end{equation*}
where $\vmu_{\parallel\vnu}$ is parallel to $\vnu$ and $\vmu_{\perp\vnu}$ is
perpendicular to $\vnu$, and $\vmu = \vmu_{\parallel\vnu} + \vmu_{\perp\vnu}$.
\end{theorem}
\begin{proof}
From Theorem \ref{geometry}, $\involute{\vmu}{\vnu}$ is the reflection of $\vmu$ in the
line defined by $\vnu$ as shown in Figure \ref{reflect}. When $\vmu$ is added to its
reflection the components of each perpendicular to $\vnu$ cancel, and the components
parallel to $\vnu$ add to give twice the stated result. The factor of ½ gives the
result as stated. Similarly, half the difference between $\vmu$ and its reflection
gives the component of $\vmu$ perpendicular to $\vnu$.
\end{proof}
Theorem \ref{projections} may be generalised to quaternions as
well as vectors. Since the scalar part of a quaternion is
invariant under an involution, the component of the quaternion
`parallel' to $\vnu$ includes the scalar part as well as the
component of the vector part parallel to $\vnu$. In other words
the `parallel' component of the quaternion is that component which
is in the same Argand plane as the axis of involution $\vnu$. The
component of the quaternion perpendicular to $\vnu$ is a vector
(the component of the vector part perpendicular to $\vnu$, and
therefore perpendicular to the Argand plane of the `parallel'
component) since the subtraction cancels out the scalar
part. As stated earlier, the
representation $a + \vmu b$ is independent of the coordinate
system in that it expresses the quaternion in terms of the
direction in 3-space of the vector part. However, the quaternion
can be rewritten in terms of a set of orthogonal basis
vectors, $\vnu_1$, $\vnu_2$, and $\vnu_3$, without recourse to a
numerical representation. The three projections across $\vnu_1$,
$\vnu_2$, and $\vnu_3$ and the conjugate anti-involution provide
the mechanism.
\newcommand{\vb}{\v{b}}
That is, we may write a quaternion $q = a + \vmu b = a + \vb$ as
\begin{equation*}
q = a + \vnu_1\alpha + \vnu_2\beta + \vnu_3\gamma = a + \vb_1 + \vb_2 + \vb_3
\end{equation*}
where $\vb_i\parallel\vnu_i$ and $\alpha$, $\beta$ and $\gamma$ are
real:
\begin{equation*}
a     = \frac{1}{2}(  q + \conjugate{q});\qquad
\vb   = \frac{1}{2}(  q - \conjugate{q});\qquad
\vb_i = \frac{1}{2}(\vb + \involute{\vb}{\vnu_i}),\quad i\in\{1,2,3\}
\end{equation*}

\end{document}